\documentclass[12pt]{article}
\usepackage[english]{babel}
\usepackage{amsfonts}

\usepackage{subcaption}
\usepackage{graphicx}
\usepackage{epstopdf}
\usepackage{geometry}
\usepackage{stmaryrd}
\usepackage{indentfirst}
\usepackage{amsthm,amsmath,amssymb}
\usepackage{mathrsfs}
\usepackage{mathtools}
\usepackage{algorithmic}
\usepackage{algorithm}
\usepackage{multirow}
\usepackage{booktabs}
\usepackage{xpatch}
\usepackage{enumitem}
\usepackage{setspace}
\usepackage{makecell}
\usepackage{tabularx}
\usepackage{adjustbox}
\usepackage{seqsplit}

\newtheorem{theorem}{\hskip 0em Theorem}[subsection]
\newtheorem{definition}{\hskip 0em Definition}[subsection]
\newtheorem{lemma}{\hskip 0em Lemma}[subsection]

\newtheorem{example}{\hskip 0em Example}[section]

\begin{document}

\title{An economic cascadic tensor multigrid method for solving high dimensional elliptic linear partial differential problems}%
\author{Jingyu Huang, Chenliang Li}%
\thanks{School of Mathematics and Computing Science, Guilin University of Electronic Technology, Guangxi Colleges and Universities Key Laboratory of Data Analysis and Computation, Center for Applied Mathematics of Guangxi (GUET)}%
\maketitle

\begin{abstract}
 In this paper, based on the tensor form, we propose a class of economic cascadic tensor multigrid method(ECTMG) for solving high dimensional elliptic linear partial differential problems. Compared with traditional methods, the new method not only reduces the storage space but also lowers the computational complexity from $\mathcal{O}(n^{3})$ to $\mathcal{O}(n^{2})$. We analyze the convergence rate of the conjugate gradient method which is based on the tensor form($\mathrm{CG}\_\mathrm{BTF}$), and then provide the convergence analysis for the new method. Finally, the effectiveness of the new method is verified through numerical examples.
\end{abstract}



\section{Introduction}
	\label{sec:introduction}
		\par Many problems in modern science and engineering involve high order partial differential equations, such as biomedical and acoustic engineering\cite{ref1}, dynamics of multi-energy systems\cite{ref2}, image restoration\cite{ref3}, mechanics of materials\cite{ref4}, and computational fluid dynamics\cite{ref5,ref6,ref7}. Partial differential equations are discretized into linear equations \eqref{eq1} by finite difference method\cite{ref8}, finite element method\cite{ref9} or spectral method\cite{ref10}.
	\begin{equation}\label{eq1}
		Ax=b,
	\end{equation}
	with $A=\left ( \sum\limits_{k=1}^{M} I_{J_{M}}\otimes \dots \otimes  I_{J_{k+1}}\otimes A^{(k)} \otimes I_{J_{k-1}}\otimes \dots \otimes I_{J_{1}}   \right), A^{(k)}\in \mathbb{R}^{J_{k}\times J_{k}}(k=1, \dots, M), A^{(k)} \otimes I_{J_{k-1}}$ represent as follows
	\begin{center}
		$A^{(k)} \otimes I_{J_{k-1}}=\begin{bmatrix}
			a_{11}I_{J_{k-1}}& \cdots  &a_{1J_{2}}I_{J_{k-1}} \\
			\vdots &  \ddots  & \vdots \\
			a_{1J_{1}}I_{J_{k-1}}& \cdots  &a_{J_{1}J_{2}}I_{J_{k-1}}
		\end{bmatrix}.$
	\end{center}
	For high-dimensional problems, as the dimension increases, the scale of the coefficient matrix $A$ increases ex\-po\-nen\-tial\-ly, which leads to the "dimension disaster". Therefore, Chen and Lu\cite{ref11} equivalently transform the linear system \eqref{eq1} into the following Sylvester tensor equation
	\begin{equation}\label{eq2}   	
		\mathcal{X}\times _{1}A^{(1)}+\mathcal{X}\times _{2}A^{(2)}+\dots +  \mathcal{X}\times _{M}A^{(M)}=\mathcal{D},  	
	\end{equation}
	with $A^{(k)}\in \mathbb{R}^{J_{k}\times J_{k}}(k=1, \dots, M) \quad$ {and} $\quad  \mathcal{D}\in \mathbb{R}^{J_{1} \times \cdots \times J_{M}}$, $\mathcal{X}\in \mathbb{R}^{J_{1}\times J_{2}\times \cdots \times J_{M}}$.  The mode-k product
	of $\mathcal{X}$ and $\mathcal{A}$ is denoted by $\mathcal{X}\times_{k}A^{(k)} $, the result is of size $J_{1}\times \cdots \times J_{k-1}\times J\times J_{k+1}\times \cdots \times J_{M}$, and its entries are
	defined by
	\begin{center}
		$(\mathcal{X}\times _{k}A^{(k)} )_{i_{1}\cdots i_{k-1}ji_{k+1}\cdots i_{M}}=\sum\limits_{i_{k}=1}^{J_{k}}x_{i_{1}\cdots i_{k-1}i_{k}i_{k+1}\cdots i_{M}}a^{(k)}_{ji_{k}} .$
	\end{center}.
	\par When $m = 2$, equation \eqref{eq2} can be equivalently transformed into the Sylvester matrix equation
	\begin{center}
		
		$A^{(1)}X+  X^{T}A^{(2)}=D,$
	\end{center}
	which is widely used in control theory\cite{ref12,ref13,ref14,ref15}, etc.
	\par Unlike equation $\eqref{eq1}$ which requires storing a large coefficient matrix $A$, equation $\eqref{eq2}$ only needs to store the smaller matrix $A^{(k)}$, significantly reducing the storage requirements. Meanwhile, the computational complexity of solving equation $\eqref{eq1}$ is $\mathcal {O}(n^{3})$, while solving equation $\eqref{eq2}$ only requires $\mathcal{O}(n^{2})$ \cite{ref16}, and the computational complexity is significantly reduced. Based on these advantages, many scholars have proposed various efficient iterative solution methods for equation $\eqref{eq2}$. A tensor format projection method and the nearest Kronecker product (NKP) preconditioner  \cite {ref11} are employed to efficiently solve the high-dimensional equation \(\eqref{eq2}\). The NKP preconditioner constructed based on the structure of the coefficient matrix can accelerate the convergence of the iterative solver, and both of these methods avoid the construction of the full coefficient matrix. Compared with the standard projection method, their floats and storage are lower, and numerical experiments have also verified their good performance. AliBeik et al.\cite{ref16} studied the tensor form of classical iterative schemes for the Sylvester tensor equation. They first derived the Arnoldi process and the full orthogonalization method (FOM) by using a product between two tensors, and then presented the tensor form of the conjugate gradient and nested conjugate gradient algorithms. A rough complexity estimation of the tensor-form Arnoldi process reveals its computational advantages over the conventional matrix-based counterpart, with numerical experiments validating the practical feasibility of the proposed tensor algorithms. Najafi-Kalyani et al.\cite{ref17} investigate global-type iterative schemes based on the Hessenberg process for solving the Sylvester tensor equation, and also study flexible variants, perturbation bounds, and Tikhonov regularization for ill-posed cases. Heyouni et al.\cite{ref18} established the tensor form of the generalized Hessenberg method for solving the Sylvester tensor equation. Huang\cite{ref19} derived the tensor form of conjugate gradient least squares method($\mathrm{CGLS}\_\mathrm{BTF}$) to solve the Sylvester tensor equation under the Tucker-product, and the numerical experiments show that the method has significant validity and superiority compared with the $\mathrm{GMRES}\_\mathrm{BTF}$ method, the $\mathrm{FOM}\_\mathrm{BTF}$ method, the$\mathrm{GLS}\_\mathrm{BTF}$ method, the GBI method and the MGBI method.
	\par The above Krylov subspace methods perform well in solving middle or small scale problems. However, they still have some limitations in large-scale problems. Recently, Chen and Li\cite{ref20} proposed a tensor multigrid method. Numerical experiments show that the proposed method exhibits extremely high computational efficiency when solving large-scale problems.
	\par  Bornemann and Deuflhard\cite{ref21} proposed a cascadic multigrid method. By simplifying the algorithm structure, this method realizes a simpler and more efficient solution strategy. Shi and Xu\cite{ref22} proposed a novel cascadic multigrid method that adopts distinct finite element spaces on the coarse and finest grids, while retaining the optimal or quasi-optimal convergence properties for second-order elliptic problems and fourth-order plate problems. Based on the classic cascadic multigrid method, Shi and Xu\cite{ref23} proposed the economic cascadic multigrid method that reduces the number of iterations at each level and improves computational efficiency. Wang and Li\cite{ref24} adopted the modulus-based matrix splitting method as a smoother and constructed a novel modulus-based cascadic multigrid(MCMG) method for elliptic variational inequality problems. The convergence of MCMG is established, and numerical experiments demonstrate its efficiency and accuracy for large-scale instances. For quasi-variational inequalities, Gao and Li\cite{ref25} introduced the modulus-based economic cascadic multigrid(MCMG) method, employing the modulus-based matrix splitting method as a smoother to accelerate convergence. A convergence analysis is provided, and numerical experiments confirm that the method achieves high efficiency and low computational cost.
	\par  Based on the advantages of the economic cascadic multigrid method in matrix form and the advantages of Sylvester tensor equation, we propose the economic cascadic tensor multigrid method(ECTMG) for solving high dimensional elliptic linear partial differential problems. Applying $\mathrm{CG}\_\mathrm{BTF}$\cite{ref16} or $\mathrm{BiCG}\_\mathrm{BTF}$\cite{ref26} as smoothers to accelerate the convergence of the iterative process. Numerical examples show that, under the same accuracy, the ECTMG method requires less storage and has a shorter runtime.
	\par  The structure of this paper is as follows. In Section 2, we provide some necessary symbols and definitions. In Section 3, we presents the economic cascadic tensor multigrid method and its convergence analysis, as well as the convergence rate analysis of the $\mathrm{CG}\_\mathrm{BTF}$ method. In Section 4, we discuss the numerical results. Finally, in Section 5, we conclude this paper.
	
	\section{Preliminaries}
\label{sec:manu}
\setcounter{subsection}{1}
\par In this paper, we use lowercase letters to denote vectors, uppercase letters to denote matrices, and Eulerian letters to denote tensors. Next, this section will briefly introduce some definitions and properties of tensors.
\begin{definition}$\textup{(\cite{ref27})}$
	let $\mathcal{A}=(a_{i_{1\times \cdots \times i_{M} }})\in \mathbb{R}^{n\times n\times \cdots\times n }  $be an $M$th-order tensor. If each of its terms $\zeta_{i_{1}i_{2}\dots i_{M}} $ satisfies
	\begin{center}
		$\zeta_{i_{1}i_{2}\dots i_{M}}=\left\{\begin{matrix}
			1,& if \quad i_{1}=\cdots =i_{M}=0,\\
			0,&otherwise,
		\end{matrix}\right.$
	\end{center}
	Then the tensor $\mathcal{A}$ is called the unit tensor, denoted by $\mathcal{I}$.
\end{definition}
\begin{definition}$\textup{(\cite{ref27})}$\label{definition2}
	For any given two tensors $\mathcal{X}, \mathcal{Y}\in \mathbb{R}^{J_{1}\times \cdots \times J_{M}}$, their inner product is defined as
	\begin{center}
		$\left \langle \mathcal{X}, \mathcal{Y}   \right \rangle=\sum\limits_{i_{1}=1}^{J_{1}}\sum\limits_{i_{2}=1}^{J_{2}} \cdots \sum\limits_{i_{m}=1}^{J_{M}}x_{i_{1}i_{2}\dots i_{M}}y_{i_{1}i_{2}\dots i_{M}}.   $
	\end{center}
	Thus, the Frobenius norm induced by the tensor inner product is expressed as
	\begin{center}
		$\left \| \mathcal{X}  \right \|=\sqrt{\left \langle \mathcal{X}, \mathcal{X}   \right \rangle }=\sqrt{\sum\limits^{J_{1}}_{i_{1}=1}\sum\limits_{i_{2}=1}^{J_{2}} \cdots \sum\limits_{i_{M}=1}^{J_{M}}x^{2}_{i_{1}i_{2}\dots i_{M}}}.$
	\end{center}
\end{definition}
\begin{definition}$\textup{(\cite{ref27})}$\label{definition3}
	The straightening operator $vec(\cdot)$ performs vectorization ($\mbox{vectorization}$) on matrices or tensors. Let $A_{mn} \in \mathbb{R}^{J_{k} \times J_{k}}$ be a matrix and $\mathcal{A} \in \mathbb{R}^{J_{1} \times \cdots \times J_{M}}$ be a tensor, then
	\begin{center}
		$vec(A)\coloneqq[a_{1}^{T}, \dots , a_{n}^{T}]^{T}, vec(\mathcal{A})\coloneqq vec(\mathcal{A}(1) ),$
	\end{center}
	where $a_{i}(i=1, 2, \dots, n)$ is the column vector of matrix $A$.
\end{definition}
\section{Main Result}
\par This section presents the economic cascadic tensor multigrid method (ECTMG).  In this method, $\mathrm{CG}\_\mathrm{BTF}$ or $\mathrm{BiCG}\_\mathrm{BTF}$ is selected as smoother. To further improve the accuracy of the solution, we use a quadratic interpolation on the $l$-th level. We also analyze the convergence rate of the $\mathrm{CG}\_\mathrm{BTF}$ algorithm and provide a convergence analysis of the ECTMG algorithm.
\subsection{The ECTMG algorithm}
\par Let the solution domain be $\Omega = \Omega_1 \times \Omega_2 \times \cdots \times \Omega_M $. Let $l$ ($l = 1, 2, \dots, L$) be the number of grid level, with $l = 1$ being the coarsest grid level. On the $l$-th level, the grid discretizes with the step size $h_l = \frac{1}{N_l + 1}$ ($N_l = 2^{l + 1} - 1$). The high-order partial differential equations are discretized by using the finite difference method, and the tensor equation system of the corresponding grid level is obtained as follows, 
\begin{equation}\label{eq3}
	\mathcal{X}\times _{1}A^{(1)}_{l}+\mathcal{X}\times _{2}A^{(2)}_{l}+\dots +  \mathcal{X}\times _{M}A^{(M)}_{l}=\mathcal{D}_{l},
\end{equation}
where, $A_{l}^{(i)} \in \mathbb{R}^{J_{k}\times J_{k}}$ ($k = 1, \dots, M$) and $\mathcal{D}_{l}$ is the right-hand term. 
\subsubsection{The number of smoothing step $m_{l}$}\label{ml}
\par In ECTMG method, the number of iteration steps is crucial for the algorithm performance. This can improve the computational efficiency. Therefore, we similar to \cite{ref23} provides the following selection rules:\\
When d=2\\
1.If $l> L_{0}$,then $m_{l}=\left [ m_{L}\beta ^{L-l} \right ]. $\\
2.If $l\le L_{0},m_{l}=[m_{*}^{\frac{1}{2}}(L-(2-\varepsilon _{0})l)\kappa _{l}].$\\
When d$\ge$3\\
1.If $l> L_{0}$,then $m_{l}=\left [ m_{L}\beta ^{L-l} \right ]. $\\
2.If $l\le L_{0},$ there are two cases:\\
(1)If $(2-\varepsilon _{0})L_{0}\le L,$then $m_{l}=[m_{*}^{\frac{1}{2}}(L-(2-\varepsilon _{0})l)\kappa _{l}].$\\
(2)If $(2-\varepsilon _{0})L_{0}>  L$,there exists a positive integer$L'< L_{0}$ such that$(2-\varepsilon _{0})L'\le   L$, for all $l\le L'$, choose $m_{l}=[m_{*}(L-(2-\varepsilon _{0})l)\kappa _{l}]$. \\	
\\ \textbf{Remark }\quad Let $\kappa_l=h_{l}^{-2}$,  and $\varepsilon _{0}$ be a positive constant in the interval [0,1].			 
\subsubsection{Prolongation}
\cite{ref28} gives the linear prolongation operator $P_{l-1}^{l}$, which is defined as follows
\begin{center}
	$P_{l-1}^{l}=\frac{1}{2}\begin{pmatrix}
		1&  &   \\
		2&  1&   \\
		1&  2&   \\
		&  1&  1\\
		&  &  \ddots 
	\end{pmatrix}_{J_{k}\times J_{k}}.$
\end{center}
Then, the tensor-form extension operator $\mathcal{P}_{l-1}^{l}$ acts on the solution $\mathcal{X}_{l-1}$ on the level $l-1$, and is defined as follows
\begin{center}
	$
	\mathcal{P}_{l-1}^{l}(\mathcal{X}_{l-1} )=  \mathcal{X}_{l-1}\times _{1}P_{l-1}^{l}\times _{2}P_{l-1}^{l} \cdots \times _{M}P_{l-1}^{l}. $
\end{center}
\par We adopted the quadratic interpolation extension operator $P_{l-1}^{l}$ as
\begin{center}
	$P_{l-1}^{l}=\begin{pmatrix}
		1 & & & & & & & \\
		\frac{3}{8} & \frac{3}{4} & -\frac{1}{8} & & & & & \\
		& & 1 & & & & & \\
		& & -\frac{1}{8} & \frac{3}{4} & \frac{3}{8} & & & \\
		& & & & 1 & & & \\
		& & & & \frac{3}{8} & \frac{3}{4} & -\frac{1}{8} & \\
		& & & & & & 1 & \\
		& & & & & & -\frac{1}{8} & \frac{3}{4} & \frac{3}{8} \\
		& & & & & & & & \ddots
	\end{pmatrix}_{(J_{k}+2)\times (J_{k}+2)}$.
\end{center}

\subsubsection{Smoother}
Let $\mathscr{T}_{l}$ be the smoother on the level $l$,and define 
\begin{center}
	$\mathcal{X}_{l}^{m_{l}}\coloneqq \mathscr{T}_{l}^{m_{l}}\mathcal{X}_{l}.$
\end{center}
\par When the coefficient matrix $A^{(k)}$ is symmetric and positive definite, we select $\mathrm{CG}\_\mathrm{BTF}$\cite{ref16} as the smoother; When the coefficient matrix $A^{(k)}$ is non-symmetric positive definite, we select $\mathrm{BiCG}\_\mathrm{BTF}$\cite{ref26} as the smoother.	
\begin{algorithm}[H]
	\caption{ECTMG algorithm}
	\label{alg:AOA}
	\renewcommand{\algorithmicrequire}{\textbf{Input:}}
	\renewcommand{\algorithmicensure}{\textbf{Output:}}
	\begin{algorithmic}[1]
		\REQUIRE Given the maximum number of mesh grid $L$, initial guess value $\mathcal{X}_{0}$, the coefficient matrices $A_{l}^{(1)}, A_{l}^{(2)}, \dots, A_{l}^{(M)}$ of the Sylvester tensor equation and the right-hand side term $\mathcal{D}_{l}$.  
		
		\STATE  On the coarsest grid level:  $\mathcal{X}_{1}$: $\mathcal{X}_{1}\coloneqq \mathscr{T}_{1}\mathcal{X}_{0}$.
		
		\FOR{$l = 2, \dots, L$}
		\STATE $\mathcal{X}_{l}\coloneqq \mathcal{P}_{l-1}^{l}\mathcal{X}_{l-1};$   
		\STATE Perform $m_{l}$ smoothing iterations on $\mathcal{X}_{l}$: $\mathcal{X}_{l}^{m_{l}}\coloneqq \mathscr{T}_{l}^{m_{l}}\mathcal{X}_{l}$;   
		\STATE Let $\mathcal{X}_{l-1} \coloneqq\mathcal{X}_{l}^{m_{l}} $
		\ENDFOR 
		\STATE  return\;$\mathcal{X}=\mathcal{X}_{L}$
	\end{algorithmic}
\end{algorithm}

\subsection{The convergence rate of $\mathrm{CG}\_\mathrm{BTF}$ method}\label{subsec:tcg-convergence}
This section presents an analysis of the convergence rate of $\mathrm{CG}\_\mathrm{BTF}$ method \cite{ref16}.
\begin{lemma}\label{definition1}
	For linear operators $\mathcal{M}:\mathbb{R}^{J_{1}\times \cdots \times J_{M}}\mapsto \mathbb{R}^{J_{1}\times \cdots \times J_{M}}$, define the mapping as follows
	\begin{center}
		$\mathcal{X}\mapsto \mathcal{M}(\mathcal{X}): = \mathcal{X}\times _{1}A^{(1)}+\dots +  \mathcal{X}\times _{M}A^{(M)}.$
	\end{center}
	Thus, equation $\eqref{eq2}$ can be rewritten as	
	\begin{center}
		$	\mathcal{M}(\mathcal{X})= \mathcal{D}.$
	\end{center} 	
	Therefore, the tensor form of the Krylov subspace\cite{ref16}
	\begin{center}
		$K_{m_{l}}(\mathcal{M}, \mathcal{R}_{0})=span\left \{\mathcal{R}_{0}, \mathcal{M}(\mathcal{R}_{0}),  \cdots, \mathcal{M}^{m_{l}-1}(\mathcal{R}_{0})\right \}$
	\end{center}
	is equivalent to the Krylov subspace 
	\begin{center}
		$K_{m_{l}}(A, r_{0})=span\left \{r_{0}, Ar_{0}, \dots, A^{m_{l}-1}r_{0}  \right \},$
	\end{center}
	where $A=\left ( \sum\limits_{k=1}^{m} I_{n^{m}}\otimes \dots \otimes  I_{n^{k+1}}\otimes A^{k} \otimes I_{n^{k-1}}\otimes \dots \otimes I_{n^{1}}   \right ) 
	, r_{0}=vec(\mathcal{R}_{0})$, $\mathcal{R}_{0}=\mathcal{D}-\mathcal{M}(\mathcal{X}_{0})$ and $\mathcal{X}_{0}\in \mathbb{R}^{J_{1}\times \cdots \times J_{M}}$ is the initial guess value of the solution to equation $\eqref{eq2}$.
	
	\par\begin{proof}
		The tensor $\mathcal{M}(\mathcal{R}_{0})$ can be vectorized as by using the operator $vec\left (  \cdot \right )$ follows:							
		\begin{equation}\label{eq4}
			vec\left (  \mathcal{M}(\mathcal{R}_{0}) \right )=\left ( \sum_{k=1}^{m} I_{n^{m}}\otimes \dots \otimes  I_{n^{k+1}}\otimes A^{k} \otimes I_{n^{k-1}}\otimes \dots \otimes I_{n^{1}}   \right ) vec(\mathcal{R}_{0}).
		\end{equation}
		
		Let $A=\left ( \sum\limits_{k=1}^{M} I_{n^{M}}\otimes \dots \otimes  I_{n^{k+1}}\otimes A^{k} \otimes I_{n^{k-1}}\otimes \dots \otimes I_{n^{1}}   \right ) 
		, r_{0}=vec(\mathcal{R}_{0}). $
		Then, \eqref{eq4} can be rewritten as
		\begin{center}
			$vec\left (  \mathcal{M}(\mathcal{R}_{0}) \right )=Ar_{0}.$	
		\end{center}
		By analogy, we can obtain
		\begin{center}
			$vec\left (  \mathcal{M}^{m_{l}-1}(\mathcal{R}_{0}) \right )=A^{m_{l}-1}r_{0}.$
		\end{center}			
		\par When $m_{l}=1$, $vec(\mathcal{R}_{0})=vec(\mathcal{M}^{0}(\mathcal{R}_{0}))=r_{0}$, so the conclusion of the lemma holds. Now assume that the conclusion of the lemma holds for $m_{l}-2$, and we will prove that it also holds for $m_{l}-1$. By using the inductive hypothesis, we have		
		\begin{align*}
			vec(\mathcal{M}^{m_{l}-1}(\mathcal{R}_{0}))
			&=vec(\mathcal{M}(\mathcal{M}^{m_{l}-2}(\mathcal{R}_{0})))\\
			&=Avec(\mathcal{M}^{m_{l}-2}(\mathcal{R}_{0}))\\
			&=A \cdot A^{m_{l}-2}vec(\mathcal{R}_{0})\\
			&=A^{m_{l}-1}vec(\mathcal{R}_{0})\\
			&=A^{m_{l}-1}r_{0}.
		\end{align*}
		
		Thus, the conclusion also holds for $m_{l}-1$.
	\end{proof}		
\end{lemma}
\par \begin{lemma}\label{lemma2}
	The linear operator $\mathcal{M}$ satisfies the following properties:
	\\$(1)\forall \mathcal{X}, \mathcal{Y}\in \mathbb{R}^{J_{1}\cdot J_{2}\cdots J_{M}\times J_{1}\cdot J_{2}\cdots J_{M} }$, $\mathcal{M}(\mathcal{X})-\mathcal{M}(\mathcal{Y})=\mathcal{M}(\mathcal{X}-\mathcal{Y}).$
	\\$(2)$For linear operators $\mathcal{M}: \mathbb{R}^{J_{1}\times \cdots \times J_{M}}\mapsto \mathbb{R}^{J_{1}\times \cdots \times J_{M}}$, $ \mathcal{M}^{2}\cdot  \mathcal{M}=\mathcal{M}\cdot \mathcal{M}^{2}$.
\end{lemma}

\begin{lemma}$\textup{(\cite{ref9,ref11})}$\label{lemma4}
	The set of all eigenvalues $\sigma (A)$ of A, which consists of the sum of all possible eigenvalues $A^{(1)}, \dots, A^{(k)}$, is expressed as
	\begin{center}
		$\sigma (A)=\sum\limits_{k=1}^{M} \sigma (A^{(k)})=\left \{ \sum\limits_{k=1}^{M}\lambda _{k}\mid \lambda _{k}\in  \sigma (A^{(k)})  \right \} $.
	\end{center}	
\end{lemma}	

\begin{lemma}$\textup{(\cite{ref16})}$\label{lemma5}
	If $A^{(k)}\in \mathbb{R}^{n^{k}\times n^{k}}(k=1, \dots, M)$ is symmetric positive definite, then the coefficient matrix A must be symmetric positive definite, but the converse is not true.
\end{lemma}
\begin{lemma}\label{lemma6}
	For $\mathcal{X}, \mathcal{Y}\in \mathbb{R}^{J_{1}\times \dots \times J_{M}} $, define the norm on $\mathbb{R}^{J_{1}\times \dots \times J_{M}}$
	\begin{center}
		$\left \langle \mathcal{X}, \mathcal{Y} \right \rangle_{\mathcal{M} }= \left \langle \mathcal{M}(\mathcal{X}), \mathcal{Y} \right \rangle$.
	\end{center}
	
\end{lemma}
Now, we present the following theorem regarding the convergence rate of $\mathrm{CG}\_\mathrm{BTF}$ method.
\begin{theorem}\label{theorem1}
	Let $A^{(k)}\in \mathbb{R}^{J_{k}\times J_{k}}(k=1, \dots, M)$ be symmetric and positive definite, with maximum and minimum eigenvalues $\lambda_{max}^{(k)}$ and $\lambda_{min}^{(k)}$, respectively, then $A \in \mathbb{R}^{J_1 \cdot J_2 \cdots J_M \times J_1 \cdot J_2 \cdots J_M}$ is symmetric and positive definite, with maximum and minimum eigenvalues $\lambda_{max}=\sum\limits_{k=1}^{M}\lambda_{max}^{(k)}$ and $\lambda_{min}=\sum\limits_{k=1}^{M}\lambda_{min}^{k}$, respectively. The approximate solution $\mathcal{X}^{(k)}$ obtained by the $\mathrm{CG}\_\mathrm{BTF}$ method satisfies:
	\begin{center}
		$\frac{\left \| \mathcal{X}^{m_{l}}-\mathcal{X}^{*} \right \|_{\mathcal{M}}^{2}  }{\left \| \mathcal{X}^{0}-\mathcal{X}^{*} \right \|_{\mathcal{M}}^{2}} \le 2\left (\frac{\sqrt{\kappa (A)}-1 }{\sqrt{\kappa (A)}+1 }   \right )^{m_{l}},$
	\end{center}
	Where $\kappa (A)=\frac{\lambda _{max}}{\lambda _{min}}$ is the (spectral) condition number of A.	
\end{theorem}
\begin{proof}
	Since $A^{(k)}\in \mathbb{R}^{J_{k}\times J_{k}}$ is symmetric and positive definite, by Lemma \ref{lemma5}, we know that $A\in \mathbb{R}^{J_{1}\cdot J_{2}\cdots J_ {M} \times J_{1}\cdot J_{2}\cdots J_{M}}$ is symmetric and positive definite. Furthermore, by Lemma \ref{lemma4}, we know that
	\begin{center}
		$\lambda _{max}=\lambda _{max}^{1} +\lambda _{max}^{2}+\dots +  \lambda _{max}^{m}$\\
		$\lambda _{min}=\lambda _{min}^{1} +\lambda _{min}^{2}+\dots +  \lambda _{min}^{m}.$
	\end{center}
	then
	$\kappa (A)=\frac{\lambda _{max}^{1} +\lambda _{max}^{2}+\dots +  \lambda _{max}^{m}}{\lambda _{min}^{1} +\lambda _{min}^{2}+\dots +  \lambda _{min}^{m}}.$\par
	Let $P_{m_{l}}$ be the set of real-coefficient polynomials of degree at most $m_{l}$. Combining Lemma \ref{definition1}, we know that
	\begin{center}
		$\mathcal{K}_{m_{l}} \left ( \mathcal{M}, \mathcal{R}_{0} \right )=\left \{ p(\mathcal{M})(\mathcal{R}_{0}), p\in P_{m_{l}-1}\right \}.$
	\end{center}
	Let $\mathcal{X}\in \mathcal{X}^{0}+\mathcal{K}_{m_{l}} \left ( \mathcal{M}, \mathcal{R}_{0} \right ), $then there exists a polynomial $p(t)\in P_{m_{l}-1}$ such that
	\begin{center}
		$\mathcal{X}=  \mathcal{X}^{0}+ p(\mathcal{M})(\mathcal{R}_{0})$.
	\end{center}
	Let $\varepsilon _{0}\overset{\bigtriangleup }{=}\mathcal{X}^{0}-\mathcal{X}^{*}, $then
	\begin{align*}
		\mathcal{X}-\mathcal{X} ^{*}&=\varepsilon _{0}+ p(\mathcal{M})(\mathcal{R}_{0})\\
		&=\varepsilon _{0}+ p(\mathcal{M})(\mathcal{D}-\mathcal{M}(\mathcal{X}^{0}))\\
		&=\varepsilon _{0}+ p(\mathcal{M})(\mathcal{M}(\mathcal{X}^{*})-\mathcal{M}(\mathcal{X}^{0}))\\
		&= \varepsilon _{0}+ p(\mathcal{M})\mathcal{M}(\mathcal{X}^{*}-\mathcal{X}^{0})&&(\mbox{lemma\ref{lemma2} (1)})\\
		&=\varepsilon _{0}+ \mathcal{M}p(\mathcal{M})(\mathcal{X}^{*}-\mathcal{X}^{0})&&(\mbox{lemma\ref{lemma2} (2)})\\
		&=(\mathcal{I}-\mathcal{M}p(\mathcal{M}))(\varepsilon _{0}).
	\end{align*}		
	According to the $vec\left (  \cdot \right )$ operator and definition \ref{definition1}, we obtain
	\begin{align*}
		vec(\mathcal{X}-\mathcal{X} ^{*})&=vec((\mathcal{I}-\mathcal{M}p(\mathcal{M}))(\varepsilon _{0}))\\
		&=vec(\varepsilon _{0})-vec(\mathcal{M}p(\mathcal{M})(\varepsilon _{0}))\\
		&=vec(\varepsilon _{0})-\sum_{i=1}^{m_{l}}A^{i}vec(\varepsilon _{0}) \\
		&=(I- \sum_{i=1}^{m_{l}}A^{i})vec(\varepsilon _{0})\\
		&=q_{m_{l}}(A)vec(\varepsilon _{0}), 
	\end{align*}
	Where $q_{m_{l}}(A)=I- \sum\limits_{i=1}^{m_{l}}A^{i}$, and $I $ is the unit matrix.
	\par Since $A$ is symmetric positive definite, we can set $A = Q\Lambda Q^{T}$, where $\Lambda = diag(\lambda _{1}, \lambda _{2}, \dots , \lambda _{n}), \lambda _{i} > 0$, and denote $y = \left [ y_{1}, y_{2}, \dots, y_{n} \right ]^{T}\triangleq Q^{T}vec(\varepsilon _{0})$.
	\begin{align*}
		\left \| \mathcal{X}^{k}-\mathcal{X}^{*}  \right \|_{\mathcal{M}}^{2}&=\underset{\mathcal{X}\in \mathcal{X}^{0}+K_{m_{l}}(\mathcal{M}, \mathcal{R}_{0}) }{min}\left \| \mathcal{X}-\mathcal{X}^{*}  \right \|_{\mathcal{M}}^{2}\\
		&=\underset{\mathcal{X}\in \mathcal{X}^{0}+K_{m_{l}}(\mathcal{M}, \mathcal{R}_{0}) }{min}\left \langle \mathcal{M}(\mathcal{X}-\mathcal{X}^{*}), \mathcal{X}-\mathcal{X}^{*}\right \rangle \\
		&=\underset{\mathcal{X}\in \mathcal{X}^{0}+K_{m_{l}}(\mathcal{M}, \mathcal{R}_{0}) }{min}(vec(\mathcal{X}-\mathcal{X}^{*}))^{T}A\;(vec(\mathcal{X}-\mathcal{X}^{*}))\\
		&=\underset{q\in \mathbb{P}_{k}, q(0)=1 }{min}vec(\varepsilon _{0})^{T}q_{m_{l}}(A)^{T}A\; vec(\varepsilon _{0})q_{m_{l}}(A)\\ 
		&=\underset{q\in \mathbb{P}_{k}, q(0)=1 }{min}\sum_{i=1}^{m_{l}}y_{i}^{2}\lambda _{i}q_{m_{l}}(\lambda _{i})^{2}\\
		&\le \underset{q\in \mathbb{P}_{k}, q(0)=1 }{min}  \underset{1\le i\le m_{l}}{max}\left \{ q_{m_{l}}(\lambda _{i})^{2} \right \} \sum_{i=1}^{m_{l}}y_{i}^{2}\lambda _{i}\\
		&=\underset{q\in \mathbb{P}_{k}, q(0)=1 }{min}\underset{1\le i\le m_{l}}{max}\left \{ q_{m_{l}}(\lambda _{i})^{2} \right \} y^{T}\Lambda y\\
		&=\underset{q\in \mathbb{P}_{k}, q(0)=1 }{min}\underset{1\le i\le m}{max}\left \{ q_{m_{l}}(\lambda _{i})^{2} \right \} vec(\varepsilon _{0})^{T}Q\Lambda Q^{T}vec(\varepsilon _{0})\\
		&=\underset{q\in \mathbb{P}_{k}, q(0)=1 }{min}\underset{1\le i\le m}{max}\left \{ q_{m_{l}}(\lambda _{i})^{2} \right \} vec(\varepsilon _{0})^{T}A\; vec(\varepsilon _{0})\\
		&=\underset{q\in \mathbb{P}_{k}, q(0)=1 }{min}\underset{1\le i\le m}{max}\left \{ q_{m_{l}}(\lambda _{i})^{2} \right \}\left \| \varepsilon _{0} \right \|^{2}_{\mathcal{M} }. 			
	\end{align*}	    
	\mbox{Then}
	\begin{center}
		$\frac{\left \| \mathcal{X}^{k}-\mathcal{X}^{*}  \right \|_{\mathcal{M}}^{2}}{\left \|\mathcal{X}^{0}-\mathcal{X}^{*}  \right \|^{2}_{\mathcal{M} }}\le \underset{q\in \mathbb{P}_{k}, q(0)=1 }{min}\underset{1\le i\le J_{1}\dots J_{M}}{max}\left \{ q_{m_{l}}(\lambda _{i})^{2} \right \}.$
	\end{center}
	\par From the properties of Chebyshev polynomials, we know that
	\begin{center}
		$\frac{\left \| \mathcal{X}^{k}-\mathcal{X}^{*}  \right \|_{\mathcal{M}}^{2}}{\left \|\mathcal{X}^{0}-\mathcal{X}^{*}  \right \|^{2}_{\mathcal{M} }}\le 2\left ( \frac{\sqrt{\kappa (A)}-1 }{\sqrt{\kappa (A)}+1}  \right )^{m_{l}}. $
	\end{center}
	Thus, we have proven the theorem.
\end{proof}
\textbf{Remark }\quad The compression factor of $\mathrm{CG}\_\mathrm{BTF}$ method is $2\left ( \frac{\sqrt{\kappa (A)}-1 }{\sqrt{\kappa (A)}+1} \right )^{m_{l}}$, which is the same as that of the traditional conjugate gradient method (CG). However, there are significant differences between the two methods. On the one hand, the  $\mathrm{CG}\_\mathrm{BTF}$ method requires less storage, which is advantageous when processing large-scale problems. On the other hand, $\mathrm{CG}\_\mathrm{BTF}$ method effectively reduces the computational complexity by adopting the Tucker product form.
\subsection{Convergence analysis}
\par Next, we present convergence analysis of Algorithm 1. First, we introduce the tensor form of the grid-related norm $\interleave \cdot  \interleave_{s,l}$.
Let $\mathcal{X}_{l} \in \mathbb{R}^{J_{1}\times \cdots \times J_{M}}   $, then
\begin{align*}
	\left \langle \mathcal{X}_{l}, \mathcal{X}_{l}  \right \rangle_{\mathcal{M}^{s}}^{2}
	&=\left \langle \mathcal{M}^{s}(\mathcal{X}_{l}), \mathcal{X}_{l}  \right \rangle^{2}\\
	&=vec(\mathcal{M}^{s}(\mathcal{X}_{l}))^{T}vec(\mathcal{X}_{l}).
\end{align*}
Similar to the proof of Theorem \ref{theorem1}, we can derive the following result
\begin{center}
	$\left \langle \mathcal{X}_{l}, \mathcal{X}_{l}  \right \rangle_{\mathcal{M}^{s}}^{2}=\sum\limits_{i=1}^{n}\lambda _{i}^{s} y_{i}^{2}  
	$,
\end{center}
where $y=\left [ y_{1}, \dots , y_{n} \right ]^{T}=Q^{T}vec(\mathcal{X}_{l} )$.
With the definition of the tensor norm, we next introduce the following definition.
\begin{definition}
	Let $\mathcal{X} \in\mathbb{R}^{J_{1}\times \cdots \times J_{M}}, $there is
	\begin{center}
		$\interleave\mathcal{X}   \interleave_{s,l} \triangleq(\sum\limits_{i=1}^{n}\lambda _{i}^{s} y_{i}^{2})^{\frac{1}{2} }=\left \langle \mathcal{X}, \mathcal{X}  \right \rangle_{\mathcal{M}_{l}^{s}}^{\frac{1}{2} }=\left \| \mathcal{M}^{\frac{s}{2} }_{l}\mathcal{X}   \right \|_{0} .$
	\end{center}
\end{definition}
According to the definition of grid-related norm, we can get the following lemma.
\begin{lemma}\label{lemma}
	Let $\mathcal{X} \in\mathbb{R}^{J_{1}\times \cdots \times J_{M}},$ there is
	\begin{center}
		$\interleave\mathcal{X}\interleave_{0,l}=\left \| \mathcal{X}  \right \|_{0} \qquad  \interleave\mathcal{X}\interleave_{1,l}=\left \langle \mathcal{X}, \mathcal{X} \right \rangle_{\mathcal{M}_{l}}^{\frac{1}{2} } $.
	\end{center}
\end{lemma}

\par Before stating the convergence theorem for the ECTMG method, we first establish a key lemma.
\begin{lemma}\label{lemma332}
	\text{If } $\mathrm{CG\_BTF}$ \text{ is} selected as the smoother, then
	\begin{center}
		$\interleave\mathscr{T}_{l}^{m_{l}}(\mathcal{X} ) \interleave_{1,l}\le \epsilon (m_{l},h_{l})\left \| \mathcal{X}   \right \|_{0} ,
		$
	\end{center}		
	where
	\begin{center}
		$\epsilon (m_{l},h_{l})=\left\{\begin{matrix}
			C\frac{h_{l}^{-1}}{m_{l}},  &m_{l}< \frac{\kappa _{l}-1}{2} , \\
			C2^{-m_{l}\frac{1}{\kappa _{l}} },&m_{l}\ge \frac{\kappa _{l}-1}{2} .
		\end{matrix}\right.$
	\end{center}	
\end{lemma}
\begin{proof}
	\begin{align*}
		\interleave\mathscr{T}_{l}^{m_{l}}(\mathcal{X} ) \interleave_{1,l}^{2}
		&=\interleave(\mathcal{I}-\mathcal{M}_{l} p(\mathcal{M}_{l}))(\mathcal{X} ) \interleave_{1,l}^{2}\\
		&=\left\langle \mathcal{M}_{l}((\mathcal{I}-\mathcal{M}_{l}p(\mathcal{M}_{l}))(\mathcal{X})),(\mathcal{I}-\mathcal{M}_{l}p(\mathcal{M}_{l})(\mathcal{X})\right\rangle\\
		&=vec(\mathcal{M}_{l}( (\mathcal{I}-\mathcal{M}_{l}p(\mathcal{M}_{l}))(\mathcal{X}))^{T}vec((\mathcal{I}-\mathcal{M}_{l}p(\mathcal{M}_{l}))(\mathcal{X}))\\
		&=vec(\mathcal{X} )^{T}(\mathcal{I}-\sum_{i=1}^{m_{l}}A_{l}^{i}  )^{T}A_{l}(\mathcal{I}-\sum_{i=1}^{m_{l}}A_{l}^{i})vec(\mathcal{X} )\\
		&=vec(\mathcal{X})^{T}q(A_{l})^{T}A_{l}q(A_{l})vec(\mathcal{X})\\
		&=y^{T}q(\Lambda )\Lambda q(\Lambda )y &&(\mbox{ let }y=Q^{T}vec(\mathcal{X}))\\
		&=\sum_{i=1}^{n}y_{i}^{2}\lambda _{i}\prod_{k=1}^{m_{l}}(1-\mu _{k}\lambda _{i})^{2}\\
		&=\sum_{i=1}^{n}y_{i}^{2}\lambda _{l,max}\cdot \frac{\lambda _{i}}{\lambda _{l,max}} \prod_{k=1}^{m_{l}}(1-\mu _{k}\lambda _{i})^{2}  \\
		&=\sum_{i=1}^{n}y_{i}^{2}\lambda _{l,max}\cdot \frac{\lambda _{i}}{\lambda _{l,max}} (1-\frac{\lambda _{i}}{\lambda _{l,max}} )^{2m_{l}} &&(\mbox{ let }\mu _{k}=\frac{1}{\lambda _{l,max}}  )\\
		&\le \lambda _{l,max}\left \{ \underset{0<x\le1}{sup}x(1-x)^{2m_{l}} \right \}\sum_{i=1}^{n}y_{i}^{2}.			      
	\end{align*}
	For $\lambda _{l,max}$($\lambda _{l,max} \mbox{ as }\mbox{the largest eigenvalue of }A_{l})$, it can be seen from the literature \cite{ref31} that
	\begin{center}
		$\lambda _{l,max}\le Ch_{l}^{-2}.$
	\end{center}
	For $\underset{{0<y\le1}}{\sup}  y(1-y)^{2m_{l}}$, it follows from the proof of Theorem \ref{theorem1} that
	\begin{center}
		$\mathcal{X}-\mathcal{X} ^{*}=(\mathcal{I}-\mathcal{M}_{l}p(\mathcal{M}_{l}))(\varepsilon _{0}),$
	\end{center}
	Let $\mathscr{T}_{l}^{m}=\mathcal{I}-\mathcal{M}_{l} p(\mathcal{M}_{l}) $,the error of the $\mathrm{CG}\_\mathrm{BTF}$ method applied to the stiffness matrix can be expressed as
	\begin{center}
		$\left \| \mathcal{X}_{j}-\mathscr{T}_{l}^{m_{l}}\mathcal{X}_{j}^{0}  \right \|_{\mathcal{M}_{l}}=\underset{q\in \mathbb{P}_{k},q(0)=1 }{min}\left \| q_{m_{l}}(\mathcal{M}_{l})(\mathcal{X}_{j}-\mathcal{X}_{j}^{0}) \right \|_{\mathcal{M}_{l} }  $,
	\end{center}
	where, $q_{m_{l}}(\mathcal{M}_{l})=\prod\limits_{k=1}^{m_{l}}(\mathcal{I}-\mu _{k}\mathcal{M}_{l}), \mu _{k}\ne 0$.
	\par Let
	\begin{center}
		$f(x)=y(1-y)^{2m_{l}}$,
	\end{center}
	where $y=\frac{\lambda _{i}}{\lambda _{l,max}} .$ When $m_{l}\ge \frac{\kappa _{l}-1}{2} $, the maximum value of $f(y)$ in the interval $(\lambda _{l,min},\lambda _{l,max})$ is
	\begin{center}
		$2^{-2m_{l}\cdot \frac{1}{\kappa _{l}} }\cdot \frac{1}{\kappa _{l}} .$
	\end{center}
	When $m_{l}< \frac{\kappa _{l}-1}{2}$, it can be seen from the literature \cite{ref23} that
	\begin{center}
		$\underset{0<y\le1}{sup}y(1-y)^{2m_{l}}\le \frac{1}{2m_{l}} .$
	\end{center}
	The lemma is proven.
\end{proof}	

\begin{lemma} \label{lemma333}
	If $m_{l}$ satisfies the selection rule in Section $\rm{\ref{ml}}$, then
	\begin{center}
		$\interleave \mathcal{X}_{L}-\mathcal{X}_{L-1}^{*}  \interleave_{1,L}\le C\frac{h_{L}}{m_{*}h_{0}}+Ch_{L}\sum\limits_{l=L_{0}+1}^{L} \frac{2^{L-l}}{\beta ^{\frac{1}{2}(L-l)m_{L}^{\frac{1}{2} } }} . $
	\end{center}
\end{lemma}

\begin{proof}
	Using the same method as in the literature \cite{ref30}, we obtain		
	\begin{align*}
		\interleave \mathcal{X}_{l}-\mathcal{X}_{l}^{*} \interleave_{1,l} &\le \interleave\mathscr{T}_{l}^{m_{l}}(\mathcal{X}_{l}-\mathcal{X}_{l-1} ) \interleave_{1,l}+ \interleave\mathscr{T}_{l}^{m_{l}}(\mathcal{X}_{l-1}-\mathcal{X}_{l-1}^{*} ) \interleave_{1,l}\\
		&= K_{1}+K_{2}.
	\end{align*}	
	For $K_{1}$, it follows from Lemma \ref{lemma332} that
	\begin{center}
		$\interleave\mathscr{T}_{l}^{m_{l}}(\mathcal{X}_{l}-\mathcal{X}_{l-1} ) \interleave_{1,l}\le \epsilon (m_{l},h_{l})\left \| \mathcal{X}_{l}-\mathcal{X}_{l-1}   \right \|_{0}.$
	\end{center}
	For $K_{2}$, it is easy to derive
	\begin{center}
		$\interleave \mathscr{T}_{l}^{m_{l}}(\mathcal{X}_{l-1}-\mathcal{X}_{l-1}^{*})  \interleave_{1,l}\le \interleave \mathcal{X}_{l-1}-\mathcal{X}_{l-1}^{*}  \interleave_{1,l}
		.$
	\end{center}
	Combining the above inequalities, we obtain
	\begin{center}
		$\interleave \mathcal{X}_{l}-\mathcal{X}_{l}^{*} \interleave_{1,l}\le \epsilon (m_{l},h_{l})\left \| \mathcal{X}_{l}-\mathcal{X}_{l-1}   \right \|_{0}+\interleave \mathcal{X}_{l-1}-\mathcal{X}_{l-1}^{*}  \interleave_{1,l}.$
	\end{center}
	Since $\mathcal{X}_{1}^{*}= \mathcal{X}_{1}$, combining Section 2.6 of Reference \cite{ref31}, we can derive through recursion that
	\begin{center}
		$\interleave \mathcal{X}_{l}-\mathcal{X}_{l}^{*} \interleave_{1,l}\le C\sum\limits _{k=2}^{l} \epsilon (m_{k},h_{k})h_{k}^{2}.$
	\end{center}
	Combining the above equation with the same proof technique used in reference \cite{ref23}, we can prove the lemma.	
\end{proof}
\subsection{Computational complexity}
For simplicity, we consider the case where $M=3$, i.e., $\mathcal{X}\in \mathbb{R}^{n\times n\times n}$. Below, we demonstrate that the computational complexity of the economic cascadic tensor multigrid (ECTMG) method is only $\mathcal{O}(n^2)$, where $n$ denotes the number of nodes in the grid. In \cite{ref16}, the computational complexity of $\mathcal{M}_{l}(X_{l})$ is $2n^{2}(\tilde{n} (A^{(1)}_{l}) + \tilde{n} (A^{(2)}_{l}) + \tilde{n} (A_{l} ^{(3)}))$, and the computational complexity of $A_{l}x_{l}$ is $2n^{3}(\tilde{n} (A_{l})),$ where $x_{l}=vec(\mathcal{X}_{l}).$ Therefore, the computational complexity of $\mathrm{CG}\_\mathrm{BTF}$ method is $2n ^{2}(\tilde {n} (A^{(1)}_{l})+\tilde{n} (A^{(2)}_{l})+\tilde{n} (A_{l}^{(3)}))$, while the computational complexity of CG method is $2n^{3}(\tilde{n} (A_{l}).$

\par According to the literature \cite{ref16}, the estimated computational complexity for the extension operator $\mathcal{P}_{l-1}^{l}(\cdot ) $ is
\begin{center}
	$\sum\limits_{i,j,k}8\tilde{n}(P_{l-1,1(i,*)}^{l} )\tilde{n}(P_{l-1,2(j,*)}^{l} )\tilde{n}(P_{l-1,3(k,*)}^{l} )\le 8n^{2}\tilde{n}(P_{l-1,1}^{l} )\tilde{n}(P_{l-1,2}^{l} )\tilde{n}(P_{l-1,3}^{l} ),   $
\end{center}
where $P_{l-1,m}^{l}$ is the linear extension operator. Therefore, the computational workloads of smoothing and prolongation are shown in Table 1.

\begin{table}[!htp]\label{jisuangongzuoliang}
	\centering
	\caption{the computational workloads of smoothing and prolongation\cite{ref21}} \label{seven}
	\begin{tabular}{lll} 
		\toprule[1.5pt] 
		& ECTMG & ECMG \\
		\midrule 
		smoother & $C_{\mathcal{S} }n^{2}$ & $C_{s} n^{3}$ \\ 
		
		prolongation & $C_{\mathcal{P} }n^{2}$ & $C_{p} n^{3}$ \\ 
		\midrule 
		Total computational workload  & $\mathcal{O}(n^{2})$ & $\mathcal{O}(n^{3})$ \\ 
		\bottomrule[1pt] 
	\end{tabular}
\end{table}
\par Combining the computational workloads of the above two processes, we present the computational workload of the ECTMG algorithm.
\begin{lemma}\label{lemma341}
	The computational cost of ECTMG method is
	\begin{center}
		$\sum\limits _{l=1}^{L}m_{l}n_{l} \le C(h_{0}^{-2}m_{*}+ (L-L_{0})m_{L} )n_{L},\quad \beta =8,$
	\end{center}
	where $n_{L}=\mathcal{O}(n^{2}).$
\end{lemma}
\begin{proof}
	\begin{equation}\label{eq5}
		\sum_{l=0}^{L} m_{l}n_{l}= \sum_{l=0}^{L_{0}} m_{l}n_{l}+\sum_{l=L_{0}+1}^{L} m_{l}n_{l}
	\end{equation}
	\par Based on the estimates of computational complexity for the smoothing process and the prolongation process discussed earlier (see Table 1), $n_{L}$ in the above equation is $n_{L}=\mathcal{O}(n^{3})$ in matrix form and $n_{L}=\mathcal{O}(n^{2})$ in tensor form.	
\end{proof}
According to Lemma \ref{lemma341}, we can see that the ECTMG algorithm requires less computational effort than the ECMG algorithm. Example 1 in Section 4 also verifies that, at the same level of accuracy, the ECTMG algorithm takes less time than the ECMG algorithm.
\par From Lemma \ref{lemma333} and Lemma \ref{lemma341}, we can obtain the following theorem.
\begin{theorem}
	The energy error of the economic cascadic tensor multigrid method is
	\begin{center}
		$\interleave \mathcal{X}_{L}-\mathcal{X}_{L}^{*} \interleave_{1,L}\le Ch_{L}(\frac{1}{m_{*}^{\frac{1}{2} }} h_{0}+C\frac{1}{m_{L}^{\frac{1}{2} }}\cdot \frac{1}{1-\frac{2}{\beta^{ \frac{1}{2} }} } ), $
	\end{center}
	And the workload is
	\begin{center}
		$\sum\limits _{l=0}^{L}m_{l}n_{l} \le C(h_{0}^{-2}m_{*}+ (L-L_{0})m_{L} )n_{L}.$
	\end{center}
\end{theorem}
	
	\section{Numerical Example}
\setcounter{section}{4}
In this section, we verify the effectiveness of the ECTMG algorithm through some numerical examples. All tests will be done with configuration: 11th Gen Intel(R) Core(TM) i5-11400H @ 2.70GHz   2.69 GHz. Let CPU(s) represent the iteration time. The number of iterations $m_{l}$ is determined by Section 3.1.1, where $\varepsilon _{0}=0.1.$ 

\par We define the error between the exact solution $ \mathcal{X}^{*} $ and the numerical solution $ \mathcal{X}^{k} $ as follows
\begin{center}
	$ E(\mathcal{X}^{*}-\mathcal{X}^{k})=\underset{i_{1},i_{2},\dots ,i_{M}}{max} \left | x_{i_{1},i_{2},\dots ,i_{M}}^{*}- x_{i_{1},i_{2},\dots ,i_{M}}^{k}\right | ,$
\end{center}
Clearly, we have $ E(\mathcal{X}^{*}-\mathcal{X}^{k})=\left \|vec(\mathcal{X}^{*}-\mathcal{X}^{k})  \right \|_{\infty } $. The symbol $\dagger$ indicates that the algorithm ran out of memory. Table 2 list the algorithms.

\begin{table}[H]
	\centering
	\caption{Algorithms}
	\label{tab:algorithms}
	\begin{tabular}{p{4cm}p{8cm}}  
		\toprule
		\textbf{Algorithms} & \textbf{Descriptions} \\
		\midrule
		ECMG(CG or BiCG) & the economic cascadic multigrid method(CG or BiCG is smoother) \\
		ECTMG($\mathrm{CG}\_\mathrm{BTF}$ or $\mathrm{BiCG}\_\mathrm{BTF}$) & 
		\makecell[t]{the economic cascadic tensor multigrid method\\ 
			($\mathrm{CG}\_\mathrm{BTF}$ or $\mathrm{BiCG}\_\mathrm{BTF}$ is smoother)} \\
		\bottomrule
	\end{tabular}
\end{table}

\begin{example}$\textup{(\cite{ref28})}$ 
	Consider the following $M$-order Poisson equation 
	\begin{center}
		$\begin{matrix}
			-\bigtriangleup u=f & in &\Omega =[0,1]^{M} \\
			u=0 & on & \partial \Omega .
		\end{matrix}
		$
	\end{center}
	\qquad Taking the partition step length as $h=\frac{1}{(n+1)} $, using the finite difference method, we can obtain the Sylvester tensor equation \eqref{eq1}, and the discrete matrix $A^{(k)}\in \mathbb{R}^{J_{k}\times J_{k}} $ is expressed as follows
	\begin{center}
		$A^{(k)}=\frac{1}{h^{2}}\begin{bmatrix}
			2&  -1&  &  & \\
			-1&  2&  -1&  & \\
			&  \ddots &\ddots   &  \ddots & \\
			&  &  -1&  2& -1\\
			&  &  &  -1&2
		\end{bmatrix}.$
	\end{center}
\end{example}

\par \textbf{ Case 1: }Set $M=3,\omega _{l}=2/\sum\limits_{k=1}^{M}(\alpha _{l}^{k}\beta _{l}^{k}), $ where $\alpha _{l}^{k}$ and $\beta _{l}^{k}$ are the maximum and minimum eigenvalues of matrix $A_{l}^{(k)}$, respectively. Select the right-hand side tensor $\mathcal{B}$ such that the solution to the corresponding 3-order Poisson equation is $u=\frac{1}{M}\cdot  {\textstyle \prod\limits_{k=1}^{M}}(x_{k}-x_{k}^{2})$, where $(x_{k})_{i}=i/(n+1),i=1,2,\dots ,n$. The table 3 lists the numerical results of the CG, $\mathrm{CG}\_\mathrm{BTF}$, ECMG(CG) and ECTMG($\mathrm{CG}\_\mathrm{BTF}$) algorithms.

\begin{table}[H]
	\centering
	\caption{Numerical results of Case 1}
	\adjustbox{max width=\textwidth}{%
		\begin{tabular}{lcccccc}
			\toprule[1.5pt]
			\multirow{2}{*}{Algorithm} 
			& \multicolumn{3}{c}{$255\times 255\times 255$}
			& \multicolumn{3}{c}{$511\times 511\times 511$} \\
			\cmidrule(lr){2-4} \cmidrule(lr){5-7}
			& CPU (s) & $E(\mathcal{X}^{*}-\mathcal{X}_{L})$ & IT 
			& CPU (s) & $E(\mathcal{X}^{*}-\mathcal{X}_{L})$ & IT \\
			\midrule
			CG       & 432.92  & 1.71e-07 & 579 & $\dagger$& $\dagger$&$\dagger$\\        
			$\mathrm{CG}\_\mathrm{BTF}$      & 186.3   & 1.71e-07 & 578 & 4223.72 & 3.84e-07 & 1187 \\        
			ECMG(CG) & 114.14   & 4.74e-07 & (4,205,333,16384,2048,256,32)
			& $\dagger$     & $\dagger$     & $\dagger$ \\    
			ECTMG($\mathrm{CG}\_\mathrm{BTF}$) & 38.80 & 4.74e-07 & 
			(3,205,256,1024,2048,256,32)& 319.75 & 3.30e-07 & (3,269,256,1024,12288,1536,192,24)\\
			\bottomrule[1.5pt]
	\end{tabular}}
\end{table}
\par As can be seen from Table 3, when n=255, the CG algorithm and the $\mathrm{CG}\_\mathrm{BTF}$ algorithm yield same error results, and the same holds for the ECMG algorithm and the ECTMG algorithm. This indicates that the CG and $\mathrm{CG}\_\mathrm{BTF}$ algorithms, as well as the ECMG and ECTMG algorithms, achieve consistent convergence performance. In terms of CPU time, the $\mathrm{CG}\_\mathrm{BTF}$ and ECTMG algorithms exhibit shorter runtime than the CG and ECMG algorithms. This is primarily because equation \eqref{eq2} in the Tucker product form reduces computational complexity, thereby decreasing the computational load and shortening the computation time. When n=511, the ECTMG and $\mathrm{CG}\_\mathrm{BTF}$ algorithms can still perform, thanks to the compact structure of the Tucker product-form of equation \eqref{eq2}, which effectively reduces memory consumption. Moreover, the ECTMG algorithm demonstrates higher computational efficiency than the $\mathrm{CG}\_\mathrm{BTF}$ algorithm in high-dimensional case.
\par\textbf{ Case 2:} Set $M=4,\omega _{l}=2/\sum\limits_{k=1}^{M}(\alpha _{l}^{k}\beta _{l}^{k}), $ where $\alpha _{l}^{k}$ and $\beta _{l}^{k}$ are the maximum and minimum eigenvalues of matrix $A_{l}^{(k)}$, respectively. Select the right-hand side tensor $\mathcal{B}$ such that the solution to the corresponding 4-order Poisson equation is $u=\frac{1}{M}\cdot  {\textstyle \prod\limits_{k=1}^{M}}(x_{k}-x_{k}^{2})$, where $(x_{k})_{i}=i/(n+1),i=1,2,\dots ,n$. The table lists the numerical results of the CG, $\mathrm{CG}\_\mathrm{BTF}$, ECMG(CG) and ECTMG($\mathrm{CG}\_\mathrm{BTF}$) algorithms.

\begin{table}[H]
	\centering
	\caption{Numerical results of Case 2}
	\adjustbox{max width=\textwidth}{%
		\begin{tabular}{lcccccc}
			\toprule[1.5pt]
			\multirow{2}{*}{Algorithm} 
			& \multicolumn{3}{c}{$63\times 63\times 63\times63$}
			& \multicolumn{3}{c}{$127\times 127\times 127\times127$} \\
			\cmidrule(lr){2-4} \cmidrule(lr){5-7}
			& CPU (s) & $E(\mathcal{X}^{*}-\mathcal{X}_{L})$ & IT 
			& CPU (s) & $E(\mathcal{X}^{*}-\mathcal{X}_{L})$ & IT \\
			\midrule
			CG       & 154.46  & 5.99e-08 & 176 & $\dagger$     & $\dagger$     & $\dagger$ \\
			
			\addlinespace[0.3em] 
			
			$\mathrm{CG}\_\mathrm{BTF}$      & 54.08   & 5.99e-08 & 175 & 6116.58          & 1.91e-07          & 334 \\
			
			\addlinespace[0.3em] 
			
			ECMG(CG) & 127.35   & 4.36e-08 & (5,77,5120,640,80)
			& $\dagger$     & $\dagger$     & $\dagger$ \\
			
			\addlinespace[0.3em] 
			
			ETCMG($\mathrm{CG}\_\mathrm{BTF}$) & 57.40 & 4.36e-08 & (4,64,256,1280,80) 
			& 2505.00       & 1.63e-07     & (4,141,256,1024,1776,111) \\
			\bottomrule[1.5pt]
	\end{tabular}}
\end{table}

\par As can be seen from Table 4, when dealing with large-scale problems, the ECTMG algorithm demonstrates its better superiority.

\begin{example} $\textup{(\cite{ref32})}$
	Consider the following three-order steady-state convection-diffusion equation
	\begin{gather*}
		\frac{\partial^{2} T}{\partial x^{2}} + \frac{\partial^{2} T}{\partial y^{2}} + \frac{\partial^{2} T}{\partial z^{2}} 
		+ \frac{\partial T}{\partial x} + \frac{\partial T}{\partial y} + \frac{\partial T}{\partial z} = 0, \\
		(x, y, z) \in \Omega = [0,1] \times [0,1] \times [0,1],			
	\end{gather*}
	and the exact solution is $T\left ( x,y,z \right )  =e^{-x}+e^{-y}+e^{-z}.$
	\par The diffusion term is discretized using the finite difference method, and the convection term is discretized using a second-order convergent scheme (Fromm scheme). The grid partition step size is set to $h=\frac{1}{n+1} $, resulting in a 3-order tensor Sylvester equation (\ref{eq2}). where the coefficient matrix $A^{(k)}\in \mathbb{R}^{J_{k}\times J_{k}}$ is as follows
	\begin{center}
		$A^{(k)} =\frac{1}{h^{2} } \begin{bmatrix}
			2&  -1&  &  & \\
			-1 &  2&  -1&  & \\
			&  \ddots &  \ddots &  \ddots & \\
			&  &  -1&  2& -1\\
			&  &  & -1 &2
		\end{bmatrix}+\frac{1}{h} \begin{bmatrix}
			3&  -5&  1&  & \\
			1 &  3&  -5&  \ddots & \\
			&  \ddots &  \ddots &  \ddots & 1\\
			&  &  1&  3& -5\\
			&  &  & 1 &3
		\end{bmatrix}.$
	\end{center}
	
	\begin{table}[H]
		\centering
		\caption{Numerical results of example 2}
		\adjustbox{max width=\textwidth}{%
			\begin{tabular}{lcccccc}
				\toprule[1.5pt]
				\multirow{2}{*}{Algorithm} 
				& \multicolumn{3}{c}{$255\times 255\times 255$}
				& \multicolumn{3}{c}{$511\times 511\times 511$} \\
				\cmidrule(lr){2-4} \cmidrule(lr){5-7}
				& CPU (s) & $E(\mathcal{X}^{*}-\mathcal{X}_{L})$ & IT 
				& CPU (s) & $E(\mathcal{X}^{*}-\mathcal{X}_{L})$ & IT \\
				\midrule
				$\mathrm{BiCG}\_\mathrm{BTF}$       & 516.38   & 3.75e-07 & 870 & 4266          & 3.84e-07      &1137 \\
				
				\addlinespace[0.3em] 
				
				ECTMG($\mathrm{BiCG}\_\mathrm{BTF}$) & 5.96 & 8.61e-08 & (9,64,256,512,64,8,1)
				& 244.41       & 1.22e-10     & (9,269,250,1024,4096,512,64,8) \\
				\bottomrule[1.5pt]
		\end{tabular}}
	\end{table}
	
\end{example}
\par As can be seen from Table 5, the ECTMG algorithm still shows favorable performance even when the coefficient matrix is non-symmetric positive definite. In comparison with the $\mathrm{BiCG}\_\mathrm{BTF}$ algorithm, the ECTMG algorithm shows distinct advantages in high-dimensional case, characterized by shorter computation time and fewer iteration steps. Moreover, the algorithm further enhances the solution accuracy by incorporating quadratic interpolation for extrapolation.

	\section{Conclusion}
	\label{sec:conclusion}
	
		In this paper, we propose the economic cascadic tensor multigrid method(ECTMG) for solving high dimensional elliptic linear partial differential problems, and analyze the convergence of the new method. This method selects $\mathrm{CG}\_\mathrm{BTF}$ or $\mathrm{BiCG}\_\mathrm{BTF}$ as smoothers to accelerate the convergence of the iterative process. Numerical examples show that, the new method requires less storage and has a shorter CPU time.\\	
	\\ \noindent {\bf Acknowledgments.} This work is supported by the Natural Science Foundation of China under Grant 12161027 and the Science and Technology Project of Guangxi (Guike AD25069086).
	

\end{document}